\input amstex
\documentstyle{amsppt}
\document

\magnification 1100

\def\a{{\alpha}}
\def\b{{\beta}}
\def\l{{\lambda}}

\def\1b{{\bold 1}}

\def\Spec{\text{Spec}\,}

\def\max{\text{max}}

\def\Gr{\text{Gr}}

\def\CC{{\Bbb C}}
\def\DD{{\Bbb D}}

\def\NN{{\Bbb N}}

\def\QQ{{\Bbb Q}}

\def\ZZ{{\Bbb Z}}

\def\and{{\quad\text{and}\quad}}

\def\qed{\hfill $\sqcap \hskip-6.5pt \sqcup$}        
\overfullrule=0pt                                    

\newdimen\Squaresize\Squaresize=14pt
\newdimen\Thickness\Thickness=0.5pt
\def\Square#1{\hbox{\vrule width\Thickness
	      \vbox to \Squaresize{\hrule height \Thickness\vss
	      \hbox to \Squaresize{\hss#1\hss}
	      \vss\hrule height\Thickness}
	      \unskip\vrule width \Thickness}
	      \kern-\Thickness}
\def\Vsquare#1{\vbox{\Square{$#1$}}\kern-\Thickness}

\title Sur l'anneau de cohomologie du sch\'ema de Hilbert de $\CC^2$\endtitle 
\rightheadtext{Anneau de cohomologie du sch\'ema de Hilbert}
\author E. Vasserot \endauthor
\abstract\nofrills{\smc R\'esum\'e.}
{Nous exprimons le produit de l'anneau de cohomologie du sch\'ema de
Hilbert en fonction du produit de l'alg\`ebre du groupe sym\'etrique. 
Cette construction est diff\'erente de celle de \cite{9, \S 4.4} 
o\`u le produit s'interpr\`ete en termes d'op\'erateurs diff\'erentiels.
Nous donnons une conjecture pour le cas des r\'esolutions cr\'epantes
de singularit\'es quotient symplectiques.}

\vskip5mm

\noindent{\smc R\'esum\'e (version anglaise).}
{We express the product of the cohomology ring of the Hilbert scheme
in terms of the center of the algebra of the symmetric group.
We give a conjecture for the case of crepant resolutions of
symplectic quotient singularities.}
\endabstract
\thanks
Recherche partiellement financ\'ee par le programme europ\'een 
no. ERB FMRX-CT97-0100.\endthanks
\endtopmatter
\document
\vskip1cm

\head A. Rappels et notations.\endhead
\subhead A.1\endsubhead
Pour tout entier positif $n$, soit $X_n$ le sch\'ema de Hilbert 
des sous-sch\'emas de longueur $n$ de $\CC^2=\Spec\CC[x,y]$. 
On pose $X_0=\{point\}$. 
Consid\'erons les courbes 
$\Sigma=\{x=0\},\Sigma'=\{y=0\}\subset\CC^2.$
Soit $Y_n\subset X_n$ le sous-ensemble des 
sch\'emas support\'es par la courbe $\Sigma.$
L'action du groupe $\CC^\times$ sur $\CC^2$ telle que
$t\cdot(x,y)=(tx,t^{-1}y)$ induit une action de 
$\CC^\times$ sur $X_n$, not\'ee $\circ$. 
Les points fixes de cette action correspondent aux id\'eaux homog\`enes 
de $\CC[x,y]$. Ils sont param\'etr\'es par l'ensemble $\Pi_n$ des
partitions de l'entier $n$.
Soit $\xi_\l\in X_n$ le point fixe associ\'e \`a la partition $\lambda$.
Soit $h(s)$ la longueur d'\'equerre 
de la case $s$ du diagramme de la partition $\l.$
Notons $h(\l)=\prod_{s\in\l}h(s)$, et notons $l(\l)$ le nombre de parts
de la partition $\l$.

\subhead A.2\endsubhead
Soit $X$ une vari\'et\'e alg\'ebrique  
munie d'une action du groupe $\CC^\times$.
Pour tout entier $k\in\NN$, soient $H^k_{\CC^\times}(X)$ et 
$H_k^{\CC^\times}(X)$ les $k$-\`emes groupes de cohomologie et d'homologie 
\'equivariantes \`a coefficients complexes. 
Nous utilisons les notations de \cite{7, \S 1}, \cite{8, \S 1} 
pour l'homologie et la cohomologie \'equivariante, i.e. 
$$H^k_{\CC^\times}(X)=H^k_{\CC^\times}(X,\CC_X),\qquad
H_k^{\CC^\times}(X)=H^{k-2\dim X}_{\CC^\times}(X,\DD_X),$$
o\`u $\DD_X$ est le dual de Verdier du faisceau constant $\CC_X$
(voir \cite{8, \S 1.13}). 
Les foncteurs $f^*,f_!$ sont d\'efinis dans \cite{7, \S 1.4}.
Les applications lin\'eaires \cite{8, 1.18.(d),(e)}
sont not\'ees respectivement
$$\cup\,:\,H_{\CC^\times}^i(X)\otimes H_{\CC^\times}^j(X)\to
H_{\CC^\times}^{i+j}(X),\qquad
\cap\,:\,H_{\CC^\times}^i(X)\otimes H^{\CC^\times}_j(X)\to
H^{\CC^\times}_{i+j}(X).$$
Rappelons que les anneaux gradu\'es
$H^*_{\CC^\times}(point)$ et $\CC[t]$ 
(o\`u $t$ est un \'el\'ement de degr\'e 2) 
sont isomorphes, ainsi que les $\CC[t]$-modules gradu\'es
$H_*^{\CC^\times}(point)$ et $\CC[t]$.
En particulier, $H^*_{\CC^\times}(X)$ est une 
$\CC[t]$-alg\`ebre gradu\'ee unitaire, 
et $H_*^{\CC^\times}(X)$ est un $\CC[t]$-module gradu\'e.
Si $X$ est une vari\'et\'e \'equidimensionelle il y a une application $\CC[t]$-lin\'eaire 
$D\,:\,H_{\CC^\times}^k(X)\to H_k^{\CC^\times}(X),$ 
qui est inversible si $X$ est lisse. 
Donc si $X$ est lisse on obtient un produit d'intersection
$$\cap\,:\,H^{\CC^\times}_i(X)\otimes H^{\CC^\times}_j(X)\to
H^{\CC^\times}_{i+j}(X).$$
Si $i\,:\,Y\hookrightarrow X$ est le plongement d'une sous-vari\'et\'e ferm\'ee 
de codimension $k$ stable par l'action du groupe $\CC^\times$,
posons $[Y]=i_!D(1_Y)\in H_{2k}^{\CC^\times}(X)$,
o\`u $1_Y\in H^0_{\CC^\times}(Y)$
est l'\'el\'ement unit\'e de l'alg\`ebre $H^*_{\CC^\times}(Y).$

\subhead A.3\endsubhead
Soit $A^{\CC^\times}_k(X)$ le $k$-\`eme groupe de Chow
\'equivariant de la $\CC^\times$-vari\'et\'e $X$ (voir \cite{3}).
L'application cycle est un morphisme de foncteurs
$cl\,:\,A_k^{\CC^\times}\to H^{\CC^\times}_{2\dim X-2k}$
(voir \cite{3, \S 2.8}, \cite{8, \S 1.17})
qui commute au produit d'intersection $\cap$, quand il est d\'efini.
A une sous-vari\'et\'e ferm\'ee $Y\subseteq X$ stable par l'action
du groupe $\CC^\times$ on associe sa classe fondamentale 
$[Y]$ dans $A^{\CC^\times}_*(X)$ (voir \cite{3, \S 2.2}). 
Par construction $cl([Y])=[Y]$. Les propri\'et\'es 
usuelles du produit d'intersection s'\'etendent 
\`a l'homologie \'equivariante. Par exemple,
si $Y_1,Y_2\subseteq X$ sont des sous-vari\'et\'es ferm\'ees
stables par l'action de $\CC^\times$ qui sont lisses et transverses au voisinage
de l'intersection $Y_1\cap Y_2$, alors $[Y_1]\cap [Y_2]=[Y_1\cap Y_2]$ dans
$H_*^{\CC^\times}(X)$.

\subhead A.4\endsubhead
Pour toute partition $\l$ soit $\l^i$ la multiplicit\'e de $i$ dans $\l$.
On pose $\l^0=\infty$.
On note $\l_1\geq\l_2\geq...$ les parts de $\l$ ordonn\'ees 
en une suite d\'ecroissante. 
Soit $\Sigma^{(\l)}$ l'adh\'erence de Zariski du sous-ensemble
$\{\sum_i\l_i[x_i]\,|\,x_i\in\Sigma,\,x_i\neq x_j\,\text{si}\,i\neq j\}$
de $X_n$. Posons
$$Y_\l=\{\xi\in X_n\,|\,\lim_{t\to\infty}t\diamond\xi\in\Sigma^{(\l)}\},\quad
X_\l=\{\xi\in X_n\,|\,\lim_{t\to 0}t\diamond\xi\in\Sigma^{(\l)}\},$$
o\`u $\diamond$ est l'action de $\CC^\times$ sur $X_n$ induite par
l'action de $\CC^\times$ sur $\CC^2$ telle que $t\diamond (x,y)=(tx,y)$. 
Notons $\theta$ le $\CC^\times$-module de dimension 1 tel que $z\in\CC^\times$
agit par multiplication par $z$.

\proclaim{Proposition A} 
(i) On a $X_n=\bigsqcup_\l X_\l$ et $Y_n=\bigsqcup_\l Y_\l$.
Les sous-ensembles $X_\l$, $Y_\l$ sont 
stables par l'action $\circ$ du groupe $\CC^\times$. 
On a $Y_\l\simeq\CC^{n}$ et $X_\l\simeq\CC^{n+l(\l)}$.
De plus $\xi_\l$ est un point r\'egulier de $X_\l$ et $Y_\l$. 
Enfin, $Y_\mu\subseteq\bar Y_\l$ si et seulement si $\l\geq\mu$.
\hfill\break
(ii) L'espace tangent $T_{\xi_\l}X_n$ est isomorphe au $\CC^\times$-module
$\bigoplus_{s\in\l}(\theta^{h(s)}\oplus\theta^{-h(s)})$, 
et $T_{\xi_\l}Y_\l$ est isomorphe \`a $\bigoplus_{s\in\l}\theta^{-h(s)}$.
\endproclaim

\noindent{\it Preuve.} Voir \cite{4} et \cite{10} par exemple.\qed

\head B. Cohomologie \'equivariante de $X_n$\endhead
\subhead B.1\endsubhead
Soit $S_n$ le groupe sym\'etrique,
et soit $S_\l\subseteq S_n$ le sous-groupe de Young 
associ\'e \`a la partition $\l$.  
Soit $Z_n$ le centre de l'alg\`ebre $\CC[S_n]$. Notons $s_\l,c_\l,h_\l\in Z_n$
le caract\`ere du module simple associ\'e \`a la partition $\l$,
la classe de conjugaison associ\'ee \`a $\l$ et le caract\`ere
du module induit $\CC[S_n]/\CC[S_\l]$. 
Soit $(m_\l;\l\in\Pi_n)$ la base duale de $(h_\l;\l\in\Pi_n)$ pour le
produit scalaire sur $Z_n$ tel que $(s_\l|s_\mu)=\delta_{\l\mu}$.
Soit $\Pi(\l,i)$ l'ensemble des partitions $\mu$ de $n+i$ telles que
$$\exists k\in[0,n]\quad\text{tel\ que}\quad
\mu^k=\l^k-1,\ \mu^{i+k}=\l^{i+k}+1,\leqno(B.2)$$
pour un certain entier $k$. Si $\mu$ et $k$ satisfont $(B.2)$, posons
$a_{\l\mu}=\mu^{i+k}.$ Soient
$\odot\,:\,Z_m\otimes Z_n\to Z_{m+n}$ et $\bullet\,:\,Z_n\otimes Z_n\to Z_n$
les applications bilin\'eaires induites par l'induction 
des repr\'esentations du groupe sym\'etrique et le produit de
l'alg\`ebre $\CC[S_n]$. On sait que 
$$c_{(i)}\odot m_\l=\sum_{\mu\in\Pi(\l,i)}a_{\l\mu}\,m_\mu,\qquad
s_\l\bullet s_\mu=\delta_{\l\mu}h(\l)s_\l.\leqno(B.3)$$ 

\subhead B.4\endsubhead
Dor\'enavant, $H^*_{\CC^\times}(X_n)$ et $H^{\CC^\times}_*(X_n)$ sont
les groupes de cohomologie et d'homologie \'equivariante relatifs \`a
l'action $\circ$.
Si $i>0$ soit $X_{n,i}$, $X'_{n,i}$ les ensembles de couples 
$(\xi,\xi')\in X_n\times X_{n+i}$
tels que $\xi\subseteq\xi'$, et le support de $\xi'-\xi$ est un point 
(de multiplicit\'e $i$) de $\Sigma$, $\Sigma'$.
Les ensembles $X_{n,i}, X'_{n,i}$ sont des sous-vari\'et\'es 
ferm\'ees de $X_n\times X_{n+i}$ 
de dimension $2n+i$ (voir \cite{5, (8.11)}).
Soit $\pi_1,\pi_2$ les projections de $X_n\times X_{n+i}$ sur chacun des 
facteurs.
Le produit de convolution par $[X_{n,i}]$ 
est l'op\'erateur $\CC[t]$-lin\'eaire
$$P_i\,:\,H^{*}_{\CC^\times}(X_n)\to\,H^{*+2i}_{\CC^\times}(X_{n+i}),
\quad\a\mapsto D^{-1}\pi_{2!}\bigl((\pi_1^*\a)\cap [X_{n,i}]\bigr)$$
(la restriction de $\pi_2$ \`a la sous-vari\'et\'e $X_{n,i}$ est propre).
Le produit de convolution par $[X_{n,i}']$, not\'e $P'_i$, est d\'efini de 
la m\^eme facon.
On a $[\Sigma']=-[\Sigma]$ dans $H^{\CC^\times}_*(\CC^2)$.
Donc $P_i'=-P_i$. 
Posons $p_\l=P_1^{\l^1}P_2^{\l^2}\cdots(1_{X_0})\in H^{2n}_{\CC^\times}(X_n)$.
Soit $\CC[t]'$ l'anneau local de $\CC[t]$ en l'id\'eal $(t-1)$.
Pour tout $\CC[t]$-module $M$, notons $M'=\CC[t]'\otimes_{\CC[t]}M$. 
D'apr\`es le th\'eor\`eme de localisation en homologie \'equivariante,
l'image directe
$${\iota_n}_!\,:\,H^{\CC^\times}_*(X_n^{\CC^\times})'\to 
H_*^{\CC^\times}(X_n)'$$
par le plongement ferm\'e $\iota_n\,:\,X_n^{\CC^\times}\hookrightarrow X_n$ 
est inversible (voir \cite{8, Proposition 4.4}).
Soit $(\,|\,)$ la forme bilin\'eaire 
$$H^*_{\CC^\times}(X_n)'\otimes_{\CC[t]'}H^*_{\CC^\times}(X_n)'\to\CC[t]',
\qquad (\a\,|\,\b)=(-1)^np_{n!}({\iota_n}_!)^{-1}D(\a\cup\b),$$ 
o\`u $p_n$ est la projection de $X_n^{\CC^\times}$ sur un point.
Soit 
$$P_i^*\,:\,H^{*}_{\CC^\times}(X_n)'\to\,H^{*-2i}_{\CC^\times}(X_{n-i})'$$
l'adjoint de l'op\'erateur $P_i$. Il est facile de voir que
$$P_i^*(\a)=D^{-1}\pi'_{1!}{(id\times\iota_n)_!}^{-1}\bigl((\pi_2^*\a)\cap [X_{n-i,i}]\bigr),$$
o\`u $\pi'_1$ est la projection $X_{n-i}\times X_n^{\CC^\times}\to X_{n-i}$.
Pour toute partition $\l\in\Pi_n$ posons 
$y_\l=D^{-1}([\bar Y_\l])\in H^{2n}_{\CC^\times}(X_n)'.$

\proclaim{Lemme 1} $(i)$ On a
$$[P_i,P_j]=[P_i^*,P_j^*]=0,\quad [P_i,P_j^*]=\delta_{i,j}i(-1)^{i-1}id.$$

$(ii)$ La famille $(y_\l; \l\in\Pi_n)$ 
est une base de $H^{2n}_{\CC^\times}(X_n)$ sur $\CC$.

$(iii)$ Il y a un isomorphisme de $\CC[t]'$-modules
$\psi\,:\,H^*_{\CC^\times}(X_n)'\to\CC[t]'\otimes Z_n,$
tel que $\psi(p_\l)=c_\l$ et $\psi(y_\l)=m_\l$. 
\endproclaim

\noindent{\it Preuve.}
Un calcul identique \`a celui de \cite{5, \S 8.4} donne
$$[P_i,P_j]=[P_i^*,P_j^*]=0,\quad [P_i,P_j^*]=\delta_{i,j}\kappa_iid,$$
o\`u $\kappa_i\in\CC$. 
D'autre part, les sous-vari\'et\'es $X_{n,i}$, $\pi_1^{-1}(\bar Y_\l)$ sont transverses
(voir la preuve du th\'eor\`eme \cite{4, Theorem 4.6}),
et l'intersection $\pi_1^{-1}(\bar Y_\l)\cap X_{n,i}$ est la r\'eunion
de sous-vari\'et\'es ferm\'ees $Z_\mu$, $\mu\in\Pi(\l,i)$, 
qui sont stables par l'action du groupe $\CC^\times$ et telles que
la restriction de $\pi_2$ \`a $Z_\mu$ est g\'en\'eriquement
un rev\^etement de $\bar Y_\mu$ de degr\'e $a_{\l\mu}$. Donc
$$P_i(y_\l)=\sum_{\mu\in\Pi(\l,i)}a_{\l\mu}\,y_\mu.\leqno(B.5)$$
Enfin, on a $(P'_i)^*([\Sigma^{(n)}])=[\Sigma^{(n-i)}]$
puisque $\Sigma$ et $\Sigma'$ sont transverses (voir \cite{5, Lemma 9.21}). 
Donc le calcul de $\kappa_i$ est le m\^eme que dans \cite{5, \S 9.3}.
La famille $(y_\l; \l\in\Pi_n)$ est une base de 
$H^*_{\CC^\times}(X_n)'$ sur $\CC[t]'$, d'apr\`es le th\'eor\`eme de localisation
et la proposition A. 
Consid\'erons l'isomorphisme de $\CC[t]'$-modules
$$\psi\,:\,H^*_{\CC^\times}(X_n)'\to\CC[t]'\otimes Z_n,
\qquad y_\l\mapsto m_\l.$$ 
Alors, $\psi(p_\l)=c_\l$ d'apr\`es $(B.5)$ et $(B.3)$. 
\qed

\subhead B.6\endsubhead
D'apr\`es la proposition A, on a 
$$H^{\CC^\times}_k(X_n)=
\bigoplus_{2n-2l(\l)+2j=k}\CC\, t^j\cap [\bar X_\l].$$
Donc 
$H_{\CC^\times}^{4n}(X_n)=t^n\cup H_{\CC^\times}^{2n}(X_n).$
L'application
$$\star\,:\,H_{\CC^\times}^{2n}(X_n)\otimes H_{\CC^\times}^{2n}(X_n)
\to H_{\CC^\times}^{2n}(X_n)$$
telle que $t^n\cup(x\star y)=x\cup y,$
munit $H_{\CC^\times}^{2n}(X_n)$ d'une structure d'alg\`ebre sur $\CC$.
Soit $[\l]\in H^{2n}_{\CC^\times}(X_n)$ l'\'el\'ement tel que
$$t^n\cup[\l]=(-1)^nh(\l)^{-1}D^{-1}([\xi_\l]).$$

\proclaim{Lemme 2} On a

$(i)$ $([\l]\,|\,[\mu])=\delta_{\l\mu},$ 

$(ii)$ $[\l]\in y_\l+\bigoplus_{\mu>\l}\CC\,y_\mu,$

$(iii)$ $(p_\l\,|\,p_\mu)=\delta_{\l\mu}z_\l$,
o\`u $z_\l=\prod_i(i^{\l^i}\l^i!)$,
\endproclaim

\noindent{\it Preuve.}
Si $X,Y$ sont des $\CC^\times$-vari\'et\'es \'equidimensionelles
telles que $Y\subseteq X$ est ferm\'ee et $X^{\CC^\times}$ est fini, alors
$$[Y]=\sum_{y\in Y^{\CC^\times}}c_y(Y)\,t^{-\dim Y}\cup[x]
\in H_*^{\CC^\times}(X)',\leqno(B.7)$$
o\`u les constantes $c_y(Y)\in\QQ$ sont telles que
si $y$ est un point r\'egulier de $Y$ alors $c_y(Y)\neq 0$, $c_y(Y)^{-1}\in\ZZ$ et 
$\theta^{1/c_y(Y)}\simeq\bigwedge^{\max}T_y^*Y$ 
(voir A.3 et \cite{2, \S 4}). 
L'\'egalit\'e $(i)$ est imm\'ediate d'apr\`es la formule de projection et la
proposition A. Plus pr\'ecis\'ement, 
si $\iota_\l$ le plongement de $\{\xi_\l\}$ dans $X_n$, alors
$$\matrix
\bigl(D^{-1}[\xi_\l]\,|\,D^{-1}[\xi_\mu]\bigr)
&=(-1)^np_{n!}{{\iota_n}_!}^{-1}\bigl([\xi_\l]\cap [\xi_\mu]\bigr)\hfill\cr
&=\delta_{\l\mu}(-1)^np_{n!}D\iota_\l^*D^{-1}\iota_{\l!}(1)\hfill\cr
&=\delta_{\l\mu}(-1)^nc_{\xi_\l}(X_n)^{-1}\,t^{2n}\hfill\cr
&=\delta_{\l\mu}t^{2n}h(\l)^2,\hfill
\endmatrix$$
o\`u la troisi\`eme \'egalit\'e est une cons\'equence de \cite{7,\S 1.10}, et 
la quatri\`eme une cons\'equence de la proposition A.
La formule $(ii)$ d\'ecoule de $(B.7)$ et de la proposition A.
L'\'egalit\'e $(iii)$ se d\'emontre comme dans \cite{4}. 
En effet, si $n\geq i$,
$$\bigl(P_i(p_\l)\,|\,p_\mu\bigr)
=\bigl(p_\l\,|\,P_i^*(p_\mu)\bigr)=i\mu^i\bigl(p_\l\,|\,p_{\mu'}\bigr)$$
d'apr\`es le lemme 1,
o\`u la partition $\mu'\in\Pi_{n-i}$ est telle que 
${\mu'}^j=\mu^j-\delta_{i,j}$ pour tout $j$.
On conclut par r\'ecurrence sur $l(\l)$.
\qed

\proclaim{Proposition B} 
$(i)$ Il existe un unique isomorphisme d'espaces vectoriels
$\phi\,:\,H^{2n}_{\CC^\times}(X_n)\to Z_n$, tel que 
$\phi([\l])=s_\l$ et $\phi(p_\l)=c_\l$.

$(ii)$ L'application $\phi$ est un morphisme d'alg\`ebres.
\endproclaim

\noindent{\it Preuve.}
D'apr\`es les lemmes 1 et 2 on a $\psi([\l])=s_\l$.
Donc, la restriction $\phi$ de $\psi$ au sous-espace
$H^{2n}_{\CC^\times}(X_n)\subset H^*_{\CC^\times}(X_n)'$  
satisfait $(i)$.
D'apr\`es la proposition A et la formule de projection on a, 
$$D^{-1}[\xi_\l]\cup D^{-1}[\xi_\mu]=
(-1)^n\delta_{\l\mu} h(\l)^2t^{2n}\cup D^{-1}[\xi_\l].$$
Donc $\phi([\l]\star[\mu])=\phi([\l])\bullet\phi([\mu])$ d'apr\`es $(B.3)$.
\qed

\head C. Th\'eor\`eme et conjecture\endhead
\subhead C.1\endsubhead
Consid\'erons les filtrations croissantes sur 
$H^{2n}_{\CC^\times}(X_n)$ et $Z_n$ telles que 
$$E_{p}=\bigoplus_{k\leq p}t^{n-k}\cup H^{2k}_{\CC^\times}(X_n),\qquad
F_{p}=\bigoplus_{n-l(\l)\leq p}\CC\, c_\l,\qquad p\in[0,n].$$
Ces filtrations sont compatibles avec $\star$ et $\bullet$ 
(voir \cite{6, \S 3.8} dans le cas de $Z_n$). 
Les anneaux gradu\'es associ\'es sont not\'es
$\Gr^{E}_*$ et $\Gr^F_*$.

\proclaim{Th\'eor\`eme} 
Les anneaux gradu\'es $H^*(X_n)$ et $\Gr_{2*}^F$ sont isomorphes.
\endproclaim

\noindent{\it Preuve.}
La famille $(p_\l; \l\in\Pi_n)$ 
est une base de $H^{2n}_{\CC^\times}(X_n)$ sur $\CC$.
D'autre part $p_\l\in t^{l(\l)}\cup H^*_{\CC^\times}(X_n)$
car $[\Sigma]=t\cup[\CC^2]$.
D'apr\`es la proposition B on a donc $\phi(E_{p})=F_{p}$. 
Les anneaux gradu\'es $\Gr^E_*$ et $\Gr^F_*$ sont donc isomorphes.
La suite spectrale de Leray 
$$E_2^{p,q}=H^p(BS^1)\otimes H^q(X_n)\Rightarrow H^{p+q}_{\CC^\times}(X_n)$$
d\'eg\'en\`ere en $E_2$, puisque les groupes cohomologie de degr\'e impair 
de $X_n$ sont nuls. 
Donc les anneaux gradu\'es $\CC[t]\otimes H^*(X_n)$ et
$$\bigoplus_{p\geq 0}t^p\cup H^*_{\CC^\times}(X_n)/
t^{p+1}\cup H^*_{\CC^\times}(X_n)$$
sont isomorphes.  En particulier on a
$\Gr^E_*\simeq H^*(X_n)$.\qed

\subhead C.2\endsubhead
Plus g\'en\'eralement, si $G\subset Sp_{2n}(\CC)$ est un sous-groupe fini
et $g\in G$, soit $w(g)$ le rang de l'op\'erateur $g-id_{\CC^{2n}}.$ 
On a $w(gh)\leq w(g)+w(h)$ pour tout $g,h$ (voir \cite{6, \S 3.8}). 
Consid\'erons la filtration croissante du centre de l'alg\`ebre $\CC[G]$
telle que $F_p$ est lin\'eairement engendr\'e par les classes de conjugaison
d'\'el\'ements $g\in G$ tels que $w(g)\leq p$.

\proclaim{Conjecture} 
Si $M\to\CC^{2n}/G$ est une r\'esolution cr\'epante,
les anneaux gradu\'es $H^*(M)$ et $\Gr_{*}^F$ sont isomorphes.
\qed
\endproclaim
 
\noindent{\bf Remarque.}
D'apr\`es \cite{1, Theorem 8.4} on sait que $\dim H^k(M)=\dim\Gr_k^F$, 
puisque $w(g)$ est le double de l'age de $g$.

\vskip2mm

\noindent{\it Remerciements.}
{\eightpoint{Je remercie V. Ginzburg qui a sugg\'er\'e l'\'enonc\'e du th\'eor\`eme.}}

\vskip1cm
{\eightpoint{
$$\matrix\format&\l&\l&\l\\
\phantom{.} & {\text{Eric Vasserot}}\\
\phantom{.} & {\text{D\'epartement de Math\'ematiques}}\\
\phantom{.} & {\text{Universit\'e de Cergy-Pontoise}}\\
\phantom{.} & {\text{2 Av. A. Chauvin}}\\
\phantom{.} & {\text{95302 Cergy-Pontoise Cedex}}\\
\phantom{.} & {\roman{France}}\\
            & {\text{email: eric.vasserot\@u-cergy.fr}}
\endmatrix$$
}}
\enddocument